\documentclass[12pt]{amsart}
\usepackage{amssymb}
\usepackage{hyperref}
\usepackage{amsmath}
\usepackage{amscd}
\usepackage{graphicx}
\usepackage{xcolor}
\usepackage{bm}
\usepackage{comment}

\DeclareMathSymbol{\shortminus}{\mathbin}{AMSa}{"39}

\newcommand{\group}[1]{\mathrm{#1}}

\newcommand{\rep}[1]{\rep{#1}}

\newcommand{\Aut}{\operatorname{Aut}}

\newcommand{\mon}{\vec{H}}

\newcommand{\M}{\group{M}}

\newcommand{\Tr}{\operatorname{Tr}}

\newcommand{\C}{\mathbb{C}}
\newcommand{\D}{\mathbb{D}}

\newcommand{\E}{\mathbb{E}}
\newcommand{\N}{\mathbb{N}}

\newcommand{\V}{\mathbf{V}}

\newcommand{\Z}{\mathbb{Z}}

\newcommand{\Y}{\mathsf{Y}}

\newtheorem{thm}{Theorem}[section]

\theoremstyle{definition}

\title[Hypergeometric Funcation of Random Matrices]{Hypergeometric Functions of Random Matrices and 
Quasimodular Forms}

\author{Jonathan Novak}

\begin{document}


\maketitle

\section{Introduction}

\subsection{Overview} 
Hypergeometric functions of complex matrices were introduced by James \cite{James}, building
on work of Bochner \cite{Bochner}, Herz \cite{Herz}, and Constantine \cite{Constantine}. These
multivariate special functions, well-known to statisticians \cite{Muirhead}, 
play many roles in random matrix theory.
The main goal of this paper is to suggest a new use for them as holomorphic observables 
of the Circular Unitary Ensemble. 
We analyze the high-dimensional behavior of the expected derivatives
of these random analytic functions, which generalize the spectral moments 
of random unitary matrices studied by Diaconis and Shahshahani \cite{DiaSha},
and show that these more general observables admit asymptotic expansions which can be described in terms of
quasimodular forms, giving an apparently new connection between the CUE and number theory.

\subsection{Orbital integral}
In this introductory section we present our main results
in an interesting special case which can be understood without too
much background. Our starting point is the orbital integral

    \begin{equation}
    \label{eqn:IntegralHCIZ}
        K_N(q,A,B) = \int_{\group{U}_N} e^{qN \Tr AVBV^*} \mathrm{d}V,
    \end{equation}

\noindent
in which the integration is over complex unitary matrices
of rank $N$ against unit mass Haar measure, with $q \in \C$ and $A,B \in \C^{N \times N}$.
This integral was studied by Harish-Chandra \cite{HC} and Itzykson-Zuber \cite{IZ}, who independently obtained 
the formula

     \begin{equation}
    \label{eqn:DeterminantalHCIZ}
        K_N(q,A,B) = \frac{1! 2! \dots (N-1)!}{(qN)^{N \choose 2}}\frac{\det[e^{qNa_ib_j}]}{\det[a_j^{N-i}] \det[b_j^{N-i}]},
    \end{equation}

\noindent
in which $a_1,\dots,a_N$ and $b_1,\dots,b_N$ are the eigenvalues of $A$ and $B$.

In Harish-Chandra's work $N$ is fixed, $q$ is imaginary, and $A,B$ are selfadjoint.
In this case, \eqref{eqn:IntegralHCIZ} is the Fourier transform of a coadjoint orbit of $\group{U}_N$, a compact 
symplectic manifold, and \eqref{eqn:DeterminantalHCIZ} can be deduced from 
exactness of stationary phase approximation \cite{DH}. In the work of Itzykson and Zuber \cite{IZ}, the objective is to estimate
\eqref{eqn:IntegralHCIZ} as $N \to \infty$, a regime in which saddle point methods do not apply
and \eqref{eqn:DeterminantalHCIZ} is not helpful.
Instead, Itzykson and Zuber opted to view $K_N$ as an entire function of 
$q,a_1,\dots,a_N,b_1,\dots,b_N$ and work with a series expansion which leverages
the symmetries of this multivariate function. 

Since $K_N$ is invariant under independent permutations of $a_1,\dots,a_N$ and
$b_1,\dots,b_N$, and also stable under swapping these two sets of
variables, its Maclaurin series can be presented in the form 

    \begin{equation}
        \label{eqn:StringExpansion}
        K_N(q,A,B) = 1+\sum_{d=1}^\infty \frac{q^d}{d!} \sum_{\alpha,\beta \in \Y^d}
       p_\alpha(A)p_\beta(B)K_N^d(\alpha,\beta),
    \end{equation}

\noindent
where the inner sum is over pairs $\alpha,\beta$ of 
Young diagrams with exactly $d$ cells, 
with $p_\alpha(A)=p_\alpha(a_1,\dots,a_N)$ and $p_\beta(B)=p_\beta(b_1,\dots,b_N)$
the corresponding Newton symmetric polynomials and $K_N^d(\alpha,\beta)$ numerical coefficients.
The same is true of $L_N= \log K_N$, which is analytic in an open neighborhood of 
the origin in $\C^{1+2N}$ and can be written

     \begin{equation}
        \label{eqn:StringExpansionConnected}
        L_N(q,A,B) = \sum_{d=1}^\infty \frac{q^d}{d!} \sum_{\alpha,\beta \in \Y^d}
        \frac{p_\alpha(A)}{(-N)^{\ell(\alpha)}}\frac{p_\beta(B)}{(-N)^{\ell(\beta)}}
        L_N^d(\alpha,\beta),
    \end{equation}

\noindent
where we have normalized each Newton polynomial by its maximum modulus on the 
unit polydisc in $\C^N$, multiplied by a sign for later notational convenience.
The numerical coefficients $K_N^d(\alpha,\beta)$ and $L_N^d(\alpha,\beta)$ contain
the same information, and are related by a version of the Exponential Formula
for series of this form \cite[Section 3]{GGN4}.
Since the asymptotic behavior of Newton polynomials is transparent, the basic problem is to analyze the $N \to \infty$ 
asymptotics of the coefficients $L_N^d(\alpha,\beta)$.

Itzykson and Zuber \cite{IZ} verified computationally that the limit 

    \begin{equation}
        L^d(\alpha,\beta) = \lim_{N \to \infty} N^{-2} L_N^d(\alpha,\beta)
    \end{equation}

\noindent
exists for all pairs of Young diagrams with $d \leq 8$ cells, and is moreover a positive integer. 
From a quantum field theory perspective, this suggests an interpretation of these limits as counting
some class of ``genus zero'' combinatorial objects, as in the case
of Hermitian matrix integrals \cite{BIPZ,ErcMcl,Guionnet:ICM}. 
Several attempts to produce such a description of $L^d(\alpha,\beta)$ led to unsatisfactorily 
contrived and complicated planar diagrams \cite{Collins,CGM,ZZ}. Eventually, it was realized
that the high-dimensional behavior of unitary matrix integrals can be naturally 
described in the language of branched coverings \cite{Novak:Banach}.

\subsection{Hurwitz numbers}
A complete $N \to \infty$ asymptotic expansion of $L_N^d(\alpha,\beta)$ is given
by the following result from \cite{GGN3}.

     \begin{thm}
    \label{thm:HCIZ}
        For each $d \in \N$ and any $\alpha,\beta \in \Y^d$, we have

            \begin{equation*}
                L_N^d(\alpha,\beta) \sim \sum_{g=0}^\infty N^{2-2g} \mon_g^d(\alpha,\beta),
            \end{equation*}

        \noindent
        as $N \to \infty$, where $\mon_g^d(\alpha,\beta)$ is 
        the number of $(r+2)$-tuples $(\pi_1,\pi_2,\tau_1,\dots,\tau_r)$ permutations from
        the symmetric group $\group{S}^d = \Aut\{1,\dots,d\}$ such that:

            \begin{enumerate}

                \smallskip
                \item 
                The product $\pi_1\pi_2\tau_1\dots\tau_r$ equals the identity
                $\iota \in \group{S}^d$;

                \smallskip
                \item 
                The group generated by $\pi_1,\pi_2,\tau_1,\dots,\tau_r$
                acts transitively on $\{1,\dots,d\}$;

                \smallskip
                \item 
                The cycle types of $\pi_1$ and $\pi_2$ are $\alpha$ and $\beta$,
                and $\tau_1,\dots,\tau_r$ are transpositions;

                \smallskip
                \item 
                We have $r=2g-2+\ell(\alpha)+\ell(\beta)$ for $g \in \N_0$;

                \smallskip
                \item 
                Writing $\tau_i = (k_i\ l_i)$ with $1 \leq k_i < l_i \leq d$, we have
                $l_1 \leq \dots \leq l_r$.
            \end{enumerate}
    \end{thm}

The coefficients $\mon_g^d(\alpha,\beta)$ in this asymptotic expansion of $L_N^d(\alpha,\beta)$
are monotone double Hurwitz numbers \cite{GGN1,GGN2}. 
The larger number $H_g^d(\alpha,\beta)$ of permutation configurations only required to satisfy
conditions $(1)-(4)$ in Theorem \ref{thm:HCIZ} are the double Hurwitz numbers studied by Okounkov
in an algebro-geometric context \cite{Okounkov:MRL}. 
Reversing the monodromy encoding, $\frac{1}{d!}H_g^d(\alpha,\beta)$ is a weighted
count of isomorphism classes of degree $d$ branched covers $f \colon \mathbf{X} \to \mathbf{Y}$
of the Riemann sphere $\mathbf{Y}$ by a compact connected Riemann surface $\mathbf{X}$ of genus $g$ having
ramification $\alpha,\beta$ over $0,\infty \in \mathbf{Y}$, simple
branching over the $r$th roots of unity, and no further singularities. 
Monotone double Hurwitz numbers also have geometric meaning
and are closely related to Euler characteristics of moduli spaces of branched coverings \cite{Novak:YM}.
However, we are concerned here with their role in an approximation problem. 
If we replace each coefficient $L_N^d(\alpha,\beta)$ in \eqref{eqn:StringExpansion} with its large $N$ asymptotic expansion
given by Theorem \ref{thm:HCIZ} and formally change order of summation, we obtain 

    \begin{equation}
    \label{eqn:TopologicalExpansion}
        L_N(q,A,B) \sim \sum_{g=0}^\infty N^{2-2g} \vec{F}_g(q,A,B)
    \end{equation}

\noindent
with

    \begin{equation}
        \vec{F}_g(q,A,B) = \sum_{d=1}^\infty \frac{q^d}{d!} \sum_{\alpha,\beta \in \Y^d}
        \frac{p_\alpha(A)}{(-N)^{\ell(\alpha)}}\frac{p_\beta(B)}{(-N)^{\ell(\beta)}}
        \mon_g^d(\alpha,\beta)
    \end{equation}

\noindent
a generating function for monotone double Hurwitz numbers of genus $g$. This 
formal expansion is satisfactory from a quantum field theory perspective,
but from the point of view of asymptotic analysis we want to know whether
\eqref{eqn:TopologicalExpansion} holds quantitatively, as an $N \to \infty$
asymptotic expansion with uniform error on a domain of the form

     \begin{equation*}
                \D_N(\varepsilon) = \{(q,A,B) \in \C^{1+2N} \colon |q|\|A\|\|B\| < \varepsilon\},
    \end{equation*}

\noindent
with $\varepsilon > 0$ an absolute constant. For this purpose it is very useful
to know a priori that the candidate orders
$\vec{F}_g(q,A,B)$ are absolutely convergent on such a domain \cite{GGN5}.

    \begin{thm}
    \label{thm:MonotoneConvergence}
        There exists a constant $\gamma > 0$ such that, for each $g \in \N_0$, 
        the multivariate power series 

            \begin{equation*}
                \vec{F}_g(q,A,B) = \sum_{d=1}^\infty \frac{q^d}{d!} \sum_{\alpha,\beta \in \Y^d}
        \frac{p_\alpha(A)}{(-N)^{\ell(\alpha)}}\frac{p_\beta(B)}{(-N)^{\ell(\beta)}}\mon_g^d(\alpha,\beta)
            \end{equation*}

        \noindent
        converges uniformly absolutely on compact subsets of the domain 

            \begin{equation*}
                \D_N(\gamma) = \{(q,A,B) \in \C^{1+2N} \colon |q|\|A\|\|B\| < \gamma\}.
            \end{equation*}
        
    \end{thm}

\noindent
Conjecturally, the optimal value of the constant in Theorem \ref{thm:MonotoneConvergence} is
$\gamma=\frac{2}{27}$. Using Theorem \ref{thm:MonotoneConvergence} together with some classical
complex function theory, one can establish the existence of $\varepsilon>0$ such that \eqref{eqn:TopologicalExpansion} holds to 
first order uniformly on compact subsets of $\D_N(\varepsilon)$ assuming only the 
necessary condition that $K_N=\exp L_N$ be nonvanishing on this domain \cite{GGN3}.

\subsection{Randomized integral}
We now consider a randomized version of $K_N(q,A,B)$ in 
which $A=U_N$ and $B=U_N^*$, where $U_N$ is a uniformly random unitary matrix of rank $N$. This yields
a random entire function of $q \in \C$, namely 

    \begin{equation}
        K_N(q,U_N,U_N^*) = \int_{\group{U}_N} e^{qN \Tr U_NVU_N^*V^*} \mathrm{d}V.
    \end{equation}

\noindent
Equivalently, by the determinantal formula \eqref{eqn:DeterminantalHCIZ} of Harish-Chandra and Itzykson-Zuber,
we have

    \begin{equation}
        K_N(q,U_N,U_N^*) = \frac{1! 2! \dots (N-1)!}{(qN)^{N \choose 2}}\frac{\det[e^{qNu_i\bar{u}_j}]}{|V(u_1,\dots,u_N)|^2},
    \end{equation}

\noindent
where

    \begin{equation}
        V(u_1,\dots,u_N) = \prod_{i<j} (u_i-u_j)
    \end{equation}

\noindent
is the Vandermonde polynomial in the eigenvalues $u_1,\dots,u_N$ of $U_N$. This shows that
$K_N(q,U_N,U_N^*)$ is sensitive to eigenvalue repulsion in the CUE, 
making the question of its $N \to \infty$ statistical behavior a natural one. 
The expected value of $K_N(q,U_N,U_N^*)$, 

    \begin{equation}
        K_N(q) = \int_{\group{U}_N \times \group{U}_N} e^{qN \Tr UVU^*V^*} \mathrm{d}(U,V)
        = 1+ \sum_{d=1}^\infty \frac{q^d}{d!} K_N^d,
    \end{equation}

\noindent
is the moment generating function of the trace of the commutator $[U_N,V_N] = U_NV_NU_N^*V_N^*$ of 
two independent uniformly random unitary matrices, scaled up by their rank. Its logarithm 

    \begin{equation}
        L_N(q) = \log K_N(q) = \sum_{d=1}^\infty \frac{q^d}{d!} L_N^d,
    \end{equation}

\noindent
which is analytic in an open neighborhood of $q=0$, is the corresponding cumulant generating
function. Recently, physicists have shown \cite{Novaes} that the limit

    \begin{equation}
        L^d = \lim_{N \to \infty} L_N^d
    \end{equation}

\noindent
exists for all $d \in \N$, and is equal to the number of commuting pairs of permutations
in the symmetric group $\group{S}^d=\mathrm{Aut}\{1,\dots,d\}$ which generate a transitive 
subgroup. These limits $L^d =\lim_N L_N$ are the counterpart of Itzykson-Zuber's limits
$L^d(\alpha,\beta)=\lim_N N^{-2} L_N^d(\alpha,\beta)$, an important difference being
that no normalization is required for $N \to \infty$ convergence in the averaged case,
and with a considerably simpler combinatorial characterization of the limit. 
The following theorem, which is a special case of the main 
result of this paper, is the counterpart of Theorem \ref{thm:HCIZ}.

    \begin{thm}
        \label{thm:MainCommutator}
       For each fixed $d \in \N$, we have

            \begin{equation*}
                L_N^d \sim \sum_{g=1}^\infty N^{2-2g} \mon_g^d
            \end{equation*}

        \noindent
        as $N \to \infty$, where $\mon_g^d$ is the number of  $(r+2)$-tuples $(\pi_1,\pi_2,\tau_1,\dots,\tau_r)$ of permutations 
        from $\group{S}^d$ such that:

               \begin{enumerate}

                \item 
                The product $\pi_1\pi_2\pi_1^{-1}\pi_2^{-1}\tau_1\dots\tau_r$ equals the identity
                $\iota \in \group{S}^d$;

                \item 
                The group generated by $\pi_1,\pi_2,\tau_1,\dots,\tau_r$
                acts transitively on $\{1,\dots,d\}$;

                \item 
                The permutations $\tau_1,\dots,\tau_r$ are transpositions;

                \item 
                We have $r=2g-2$ for $g \in \N$;

                \item 
                Writing $\tau_i = (k_i\ l_i)$ with $1 \leq k_i < l_i \leq d$, we have
                $l_1 \leq \dots \leq l_r$.
            \end{enumerate}

        \end{thm}

The coefficients $\mon_g^d$ in this asymptotic expansion of $L_N^d$ 
are monotone simple Hurwitz numbers with torus target. The larger number $H_g^d$ of permutation configurations only 
required to satisfy conditions $(1)-(4)$ in Theorem \ref{thm:MainCommutator} is again a classical Hurwitz number: 
$\frac{1}{d!}H_g^d$ is a weighted count of isomorphism classes of degree $d$ branched covers
$f \colon \mathbf{X} \to \mathbf{Y}$ of an elliptic curve $\mathbf{Y}$ by a compact connected Riemann 
surface $\mathbf{X}$ of genus $g$ with $2g-2$ simple branch points at specified locations. 
Thus, taking a disordered version of the orbital integral \eqref{eqn:IntegralHCIZ} and averaging
the disorder, what emerges as $N \to \infty$ is the formal topological expansion

    \begin{equation}
    \label{eqn:EllipticTopologicalExpansion}
        L_N(q) = \log \int_{\group{U}_N \times \group{U}_N} e^{qN \Tr UVU^*V^*} \mathrm{d}(U,V) 
        \sim \sum_{g=1}^\infty N^{2-2g} \vec{F}_g(q)
    \end{equation}

\noindent
in which

    \begin{equation}
    \label{eqn:MonotoneGenusSpecific}
        \vec{F}_g(q) = \sum_{d=1}^\infty \frac{q^d}{d!} \mon_g^d
    \end{equation}

\noindent
is the generating function for monotone simple Hurwitz numbers of genus $g$
with elliptic as opposed to rational base curve.
This is the counterpart of \eqref{eqn:TopologicalExpansion}, but now with an $O(1)$ leading term
given explicitly by

    \begin{equation}
        \vec{F}_1(q) = \sum_{n=1}^\infty \log \frac{1}{1-q^n}.
    \end{equation}

Once again, from an analytic perspective we would like to know that the 
formal topological expansion \eqref{eqn:EllipticTopologicalExpansion} 
holds quantitatively, as an $N \to \infty$ asymptotic expansion of the 
univariate analytic function with uniform error for $|q| \leq \varepsilon$
with $\varepsilon > 0$ an absolute constant. It is again very useful to 
know in advance that the candidate orders $\vec{F}_g(q)$ in this expansion 
are analytic on a common neighborhood of $q=0$. This is confirmed by 
the following counterpart of Theorem \ref{thm:MonotoneConvergence}, which is a special case of 
a more general result established below.

    \begin{thm}
        \label{thm:QuasimodularCommutator}
        For each $g \geq 2$, the generating function $\vec{F}_g(q)$ is a polynomial in 

            \begin{equation*}
            \begin{split}
              E_2(q) &= 1 -24\sum_{n=1}^\infty \frac{nq^n}{1-q^n}, \\
              E_4(q) &= 1 +240\sum_{n=1}^\infty \frac{n^3q^n}{1-q^n}, \\
               E_6(q) &= 1 -504\sum_{n=1}^\infty \frac{n^5q^n}{1-q^n}.
              \end{split}
            \end{equation*}

    \end{thm}

It is a very classical fact that the Eisenstein series

    \begin{equation}
        E_{2k}(q) = 1 + \frac{2}{\zeta(1-2k)}\sum_{n=1}^{\infty} \frac{n^{2k-1}q^n}{1-q^n}
    \end{equation}

\noindent
is the Fourier expansion of a modular form for $\group{SL}_2(\Z)$ provided $k \geq 2$. However,
$E_2(q)$ is not modular. The algebra of functions on the open unit disc generated by all $E_{2k}(q)$
including the case $k=1$ was studied by Kaneko and Zagier \cite{KanZag}, who called
its elements quasimodular forms and showed that $E_2,E_4$, and $E_6$ are 
independent generators. Thus, Theorem \ref{thm:QuasimodularCommutator} says that the 
generating function $\vec{F}_g(q)$ for monotone simple Hurwitz numbers with torus
target is a quasimodular form for $g \geq 2$. 

Theorem \ref{thm:QuasimodularCommutator} is the 
``monotone analogue'' of a famous theorem of Dijkgraaf \cite{Dijk}, who showed
that Theorem \ref{thm:QuasimodularCommutator} holds verbatim for

    \begin{equation}
        F_g(q) = \sum_{d=1}^\infty \frac{q^d}{d!} H_g^d.
    \end{equation}

\noindent
In fact, Dijkgraaf also arrived at this result from a large $N$ perspective. 
It was shown by Gross and Taylor \cite{GrossTayl} that the chiral free 
energy of $\group{U}_N$ Yang-Mills theory with spacetime a torus of area
$t>0$ admits the formal $N \to \infty$ expansion 

    \begin{equation}
        \log Z_N^+(t) \sim \sum_{g=1}^\infty \left( \frac{N}{t} \right)^{2-2g}
        F_g(q), \quad q= e^{-\frac{t}{2}},
    \end{equation}

\noindent
and Dijkgraaf's result gives a large $N$ mirror symmetry to accompany this
large $N$ gauge-string duality for 2D Yang-Mills theory in the Calabi-Yau case.

\section{Hypergeometric Functions of the CUE}
In this paper we show that Theorems \ref{thm:MainCommutator} and 
\ref{thm:QuasimodularCommutator} hold in much greater generality, namely 
for all hypergeometric functions of the Circular Unitary Ensemble.
We construct these random analytic functions in this section.

    \subsection{Univariate hypergeometric functions}
    Let

    \begin{equation}
    \label{eqn:UniveriateHypergeometric}
        G\left(z\, \bigg{|} \begin{matrix} x_1 & \dots & x_m \\
        y_1 & \dots & y_n \end{matrix}\right) = 1+\sum_{d=1}^\infty z^d 
        \frac{(x_1)_d \dots (x_m)_d}{(y_1)_d \dots (y_n)_d},
    \end{equation}

\noindent
be the general hypergeometric series,
in which $(x)_d=x(x+1) \dots (x+d-1)$ is the Pochammer polynomial of degree 
$d$. This is a well-defined formal power series in $z$
provided the lower parameters $y_1,\dots,y_n \in \C$ are barred from $\{0,-1,-2,\dots\}$
to avoid division by zero in its coefficients,

    \begin{equation}
         \frac{(x_1)_d \dots (x_m)_d}{(y_1)_d \dots (y_n)_d} = 
         \prod_{c=0}^{d-1} \frac{(x_1+c) \dots (x_m+c)}{(y_1+c) \dots (y_n+c)}.
    \end{equation}

\noindent
If no parameters are present, then \eqref{eqn:UniveriateHypergeometric} is 
the ordinary geometric series. 
If any of the upper parameters $x_1,\dots,x_m \in \C$ is a nonpositive integer
the hypergeometric series is a polynomial. Otherwise, 
it has radius of convergence
equal to zero if $m>n$, one if $m=n$, and infinity if $m <n$. 
Hypergeometric functions are those analytic functions of a complex 
variable $z$ whose Maclaurin series is a
hypergeometric series. Basic examples are the Gaussian hypergeometric function,

    \begin{equation}
         G\left(z\, \bigg{|} \begin{matrix} x_1 & x_2 \\
        y & 1 \end{matrix}\right) = 1 + \sum_{d=1}^\infty z^d
        \frac{(x_1)_d(x_2)_d}{(y)_d(1)_d},
    \end{equation}

\noindent
the binomial series

    \begin{equation}
        G\left(z\, \bigg{|} \begin{matrix} x \\
        1 \end{matrix}\right) = 1 + \sum_{d=1}^\infty z^d \frac{(x)_d}{(1)_d},
    \end{equation}

\noindent
and the exponential series

    \begin{equation}
        G\left(z\, \bigg{|} \begin{matrix} {} \\
        1 \end{matrix}\right) = 1 + \sum_{d=1}^\infty z^d \frac{1}{(1)_d}.
    \end{equation}

\subsection{Multivariate hypergeometric functions}
Hypergeometric functions are ubiquitous and
the question of multivariate generalizations is a natural one.
Indeed, there are several distinct constructions of multivariate
hypergeometric functions, some quite sophisticated. We are interested in the
approach taken by James in the statistics literature \cite{James}, which has a very classical flavor
and turns out to be useful in applications.

The basic ingredient in James's construction is Cauchy's ``bialternant'' 

    \begin{equation}
    \label{eqn:Bialternants}
        s_\lambda(x_1,\dots,x_N) = \frac{\det [x_i^{N-j+\lambda_j}]}{\det [x_i^{N-j}]},
    \end{equation}

\noindent
a symmetric polynomial in $N$ variables associated
to every Young diagram $\lambda$ with at most $N$ rows.
The bialternants form a linear basis of the algebra of symmetric polynomials in $N$ variables
as $\lambda$ ranges over the set $\mathrm{Y}_N$ of all Young diagrams with at most
$N$ rows. Schur demonstrated that this basis consists precisely of
the characters of the irreducible homogeneous polynomial representations
of $\mathrm{GL}_N=\mathrm{Aut} \C^N$, and consequently bialternants  
are today known as Schur polynomials \cite{Macdonald}. 

Cauchy also considered the multivariate power series

    \begin{equation}
    \label{eqn:MultivariateGeometric}
        G_N(z,A,B) = 1 + \sum_{d=1}^\infty z^d \sum_{\lambda \in \mathrm{Y}_N^d} s_\lambda(A)s_\lambda(B),
    \end{equation}

\noindent
where the internal sum is over Young diagrams with at most $N$ rows and exactly 
$d$ cells, and $s_\lambda(A)=s_\lambda(a_1,\dots,a_N)$ and $s_\lambda(B)=s_\lambda(b_1,\dots,b_N)$ 
are two copies of the same Schur polynomial in algebraically independent sets of 
commuting variables. He showed that 

    \begin{equation}
    \label{eqn:CauchyIdentity}
        G_N(z,A,B) = \prod_{i,j=1}^N \frac{1}{1-za_ib_j},
    \end{equation}

\noindent
identifying $G_N$ as a nontrivial multivariate extension of the geometric series
which converges absolutely on

    \begin{equation}
        \D_N = \{(z,A,B) \in \C^{1+2N} \colon |z|\|A\|\|B\| < 1\},
    \end{equation}

\noindent
where $\|A\| = \max\{|a_1|,\dots,|a_N|\}$ denotes the $\ell^\infty$-norm of 
$A \in \C^N$. Cauchy's identity \eqref{eqn:CauchyIdentity} can be proved in 
several different ways; see \cite{Bump,Macdonald,Stanley}.

A theory of multivariate hypergeometric series built on the Cauchy identity 
must be a theory of multivariate functions of the form

    \begin{equation}
    \label{eqn:MultivariateHypergeometric}
        G_N\left(z,A,B\, \bigg{|}\, \begin{matrix} x_1 & \dots & x_m \\
        y_1 & \dots & y_n \end{matrix}\right) = 1+\sum_{d=1}^\infty z^d \sum_{\lambda \in 
        \mathrm{Y}_N^d}
        s_\lambda(A)s_\lambda(B)\frac{(x_1)_\lambda \dots (x_m)_\lambda}{(y_1)_\lambda \dots (y_n)_\lambda},
    \end{equation}

\noindent
where $(x)_\lambda$ is an extension of the Pochammer polynomial from integers 
to integer partitions which remains to be defined. A natural choice is 

    \begin{equation}
        (x)_\lambda = \prod_{\Box \in \lambda} (x+c(\Box)),
    \end{equation}

\noindent
the content polynomial of $\lambda$, 
which plays a basic role in symmetric function theory \cite[Chapter 1, Example 11]{Macdonald}.
Here $c(\Box)$ is the column index minus the row index of a given cell $\Box \in \lambda$, 
so that $(x)_\lambda$ reverts to $(x)_d$ if $\lambda$ is a row of $d$ cells. 

We take \eqref{eqn:MultivariateHypergeometric} as our
definition of the general multivariate hypergeometric series. 
When no parameters are present, \eqref{eqn:MultivariateHypergeometric}
is Cauchy's multivariate generalization of the geometric series.
As a formal power series in $1+2N$ commuting indeterminates, \eqref{eqn:MultivariateHypergeometric} is well-defined
provided the lower parameters $y_1,\dots,y_n \in \C$ are barred from 
$\{N-1,N-2,N-3,\dots\}$ to avoid division by zero in its coefficients,

    \begin{equation}
         \frac{(x_1)_\lambda \dots (x_m)_\lambda}{(y_1)_\lambda \dots (y_n)_\lambda} = 
         \prod_{\Box \in \lambda} \frac{(x_1+c(\Box)) \dots (x_m+c(\Box))}{(y_1+c(\Box)) \dots (y_n+c(\Box))}.
    \end{equation}

\noindent
If any of the upper parameters $x_1,\dots,x_n \in \C$ is an integer less than $N$ 
the multivariate hypergeometric series is a multivariate polynomial.
Otherwise, it does not converge absolutely on any open neighboorhood of the 
origin in $\C^{1+2N}$ if $m>n$, converges absolutely on $\D_N$ but not on 
any proper open superset thereof if $m=n$, and converges absolutely on 
all of $\C^{1+2N}$ if $m<n$. We define multivariate hypergeometric functions 
to be those analytic functions of $1+2N$ complex variables $z,a_1,\dots,a_N,b_1,\dots,b_N$ whose 
Maclaurin series is a multivariate hypergeometric series.

\subsection{Hypergeometric functions of matrices}
Since Schur polynomials
are irreducible characters of the general linear group, it is very 
natural to think of $A$ and $B$ as a pair of $N \times N$ complex matrices, 
with the understanding that \eqref{eqn:MultivariateHypergeometric} depends only on the
eigenvalues $a_1,\dots,a_N$ and $b_1,\dots,b_N$ of these matrices. 
This perspective is also natural in multivariate statistics and random matrix theory, 
where the multivariate hypergeometric functions we are considering
are usually called ``hypergeometric functions of two complex matrix arguments''
\cite{GrossRich}. This is an awkward name but a useful point of view: for example, 

    \begin{equation}
        G_N(z,A,B) = \int_{\mathrm{U}_N} \frac{\mathrm{d}V}{\det(I-zAVBV^*)}
    \end{equation}

\noindent
for any $A,B \in \C^{N \times N}$.
This integral representation 
was certainly not known to Cauchy and is apparently due to James \cite{James}, 
who showed more generally that

    \begin{equation}
    \label{eqn:MultivariateBinomial}
        G_N\left(z,A,B\, \bigg{|}\, \begin{matrix} x \\
        N \end{matrix}\right) 
        = \int_{\mathrm{U}_N} \frac{\mathrm{d}V}{\det(I-zAVBV^*)^x}, 
    \end{equation}

\noindent
identifying this hypergeometric function of two matrices as the natural multivariate generalization 
of the binomial series. This is already an interesting result 
with many applications. For example, in a recent paper of Marcus, Spielman,
and Srivastava \cite{MSS} one finds the integral evaluation 

    \begin{equation}
        \int_{\group{U}_N} \det(I-zAUBU^*) \mathrm{d}U = 1 + \sum_{d=1}^N 
        z^d \frac{e_d(A)e_d(B)}{{N \choose d}},
    \end{equation}

\noindent
where $e_d(A)=e_d(a_1,\dots,a_N)$ and $e_d(B) = e_d(b_1,\dots,b_N)$ are elementary 
symmetric polynomials of degree $d$ in the eigenvalues of $A,B \in \C^{N \times N}$. It is 
easy to see that this is just the case $x=-1$ of \eqref{eqn:MultivariateBinomial}.

If 

    \begin{equation}
         G_N\left(z,A,B\, \bigg{|}\, \begin{matrix} x \\
        N \end{matrix}\right) = 1 + \sum_{d=1}^\infty z^d 
        \sum_{\lambda \in \Y_N^d} s_\lambda(A) s_\lambda(B) \frac{(x)_\lambda}{(N)_\lambda}
    \end{equation}

\noindent
is the multivariate binomial, then 

    \begin{equation}
    \label{eqn:SerialHCIZ}
         G_N\left(z,A,B\, \bigg{|}\, \begin{matrix} {} \\ N \end{matrix}\right)
         = 1+\sum_{d=1}^\infty z^d \sum_{\lambda \in \Y_N^d}
        s_\lambda(A)s_\lambda(B) \frac{1}{(N)_\lambda}
    \end{equation}

\noindent
must be the multivariate exponential, and indeed

    \begin{equation}
         G_N\left(z,A,B\, \bigg{|}\, \begin{matrix} {} \\ N \end{matrix}\right) =
         \int_{\group{U}_N} e^{z\Tr AUBU^*} \mathrm{d}U,
    \end{equation}

\noindent
which connects the multivariate hypergeometric functions defined in this 
section to the material in the previous section.
Another important example is the multivariate Bessel function 

    \begin{equation}
        G_N\left(z,A,B\, \bigg{|}\, \begin{matrix} {} & {} \\M & N \end{matrix}\right)
         = 1+\sum_{d=1}^\infty z^d \sum_{\lambda \in \Y_N^d}
        s_\lambda(A)s_\lambda(B) \frac{1}{(M)_\lambda(N)_\lambda}
    \end{equation}

\noindent
which, when $M=N=1$, reduces to $I_0(2\sqrt{zab})$, where 

    \begin{equation}
        I_0(z) = 1 + \sum_{d=1}^\infty \left(\frac{z}{2}\right)^2 \frac{1}{d!d!}
        = \int_{\group{U}_1} e^{zu+z\bar{u}} \mathrm{d}u
    \end{equation}

\noindent
is the classical Bessel function which, when $z=i$, defines the Hankel transform. 
The multivariate Bessel function has the matrix integral representation \cite{Novak:BK}

    \begin{equation}
        G_N\left(z,A,B\, \bigg{|}\, \begin{matrix} {} & {} \\M & N \end{matrix}\right)
        = \int_{\group{U}_M \times \group{U}_N} e^{z_1\Tr A_1VB_1U^* + z_2\Tr A_2UBV^*} 
        \mathrm{d}(U,V),
    \end{equation}

\noindent
where $z_1,z_2 \in \C$ are any complex numbers with $z_1z_2=z$ and 
$A_1,A_2,B_1,B_2 \in \C^{M \times N}$ are any complex matrices with $A_1^*A_2=A$ and 
$B_1^*B_2=B$. According to a recent result of McSwiggen and the author \cite{McNovak}
every multivariate hypergeometric function whose lower parameters
$y_1,\dots,y_n$ are integers can be represented as an $n$-fold unitary 
matrix integral.

\subsection{Hypergeometric functions of the CUE}
It is also true that every multivariate
hypergeometric function can be represented as a ratio of determinants.
This is a remarkable result of
Khatri \cite{Khatri} and Gross-Richards \cite{GrossRich}, who showed
that the general multivariate hypergeometric series \eqref{eqn:MultivariateHypergeometric} is,
up to simple explicit factors, equal to 

    \begin{equation}
    \label{eqn:KhatGrossRich}
          \frac{\det\left[G\left(za_ib_j\, \bigg{|} \begin{matrix} x_1-N+1 & \dots & x_m-N+1 \\
        y_1-N+1 & \dots & y_n-N+1 \end{matrix}\right)\right]}{\det[a_j^{N-i}] \det[b_j^{N-i}]}.
    \end{equation}

\noindent
This immediately implies the convergence conditions on 
the general multivariate hypergeometric series stated above. Moreover, 
the determinantal representation shows that every hypergeometric
function of two matrix arguments gives a holomorphic observable
of the CUE: the random univariate function 

    \begin{equation}
    \label{eqn:GeneralObservable}
         G_N\left(z,U_N,U_N^*\, \bigg{|}\, \begin{matrix} x_1 & \dots & x_m \\
        y_1 & \dots & y_n \end{matrix}\right)
    \end{equation}

\noindent
is proportional to

    \begin{equation}
    \label{eqn:RandomDeterminant}
         \frac{\det\left[G\left(zu_i\bar{u}_j\, \bigg{|} \begin{matrix} x_1-N+1 & \dots & x_m-N+1 \\
        y_1-N+1 & \dots & y_n-N+1 \end{matrix}\right)\right]}{|V(u_1,\dots,u_N)|^2},
    \end{equation}

\noindent
with eigenvalue repulsion visible in the denominator. The expected value

    \begin{equation}
        G_N\left(z\, \bigg{|}\, \begin{matrix} x_1 & \dots & x_m \\
        y_1 & \dots & y_n \end{matrix}\right) = \mathbb{E}G_N\left(z,U_N,U_N^*\, \bigg{|}\, \begin{matrix} x_1 & \dots & x_m \\
        y_1 & \dots & y_n \end{matrix}\right)
    \end{equation}

\noindent
of this random analytic function is a deterministic analytic function whose
derivatives at $z=0$ are analogous to moments. Indeed, these derivatives are linear combinations
of the expectations 

    \begin{equation}
        \label{eqn:Diaconis}
        \E[p_\alpha(U_N) p_\beta(U_N^*)] = \int_{\group{U}_N} p_\alpha(U) p_\beta(U^*) \mathrm{d}U
    \end{equation}

\noindent
studied by Diaconis-Shahshahani \cite{DiaSha} and Diaconis-Evans \cite{DiaEva}, 
who showed that

    \begin{equation}
        \label{eqn:DiaconisLimit}
        \lim_{N \to \infty} \E_\alpha(U_N) p_\beta(U_N^*)] = \delta_{\alpha\beta} \frac{d!}{|C_\alpha|}
    \end{equation}

\noindent
where $C_\alpha$ is the conjugacy class of permutations of cycle type $\alpha$ in 
$\group{S}^d$. This result is equivalent to the Cauchy identity, and in this 
paper we generalize it to arbitrary hypergeometric functions of the CUE as defined above.

\section{Asymptotics of Expected Derivatives}
In this section we give the general version of Theorem \ref{thm:MainCommutator}.

    \subsection{Scaling considerations}
    Let us return to the general multivariate hypergeometric series
    
        \begin{equation}
        G_N\left(z,A,B\, \bigg{|}\, \begin{matrix} x_1 & \dots & x_m \\
        y_1 & \dots & y_n \end{matrix}\right) = 1+\sum_{d=1}^\infty z^d \sum_{\lambda \in 
        \mathrm{Y}_N^d}
        s_\lambda(A)s_\lambda(B)\frac{(x_1)_\lambda \dots (x_m)_\lambda}{(y_1)_\lambda \dots (y_n)_\lambda}
        \end{equation}

    \noindent
    and scale it for large $N$. Let
    $u_1,\dots,u_m,v_1,\dots,v_n \in \C$ be complex parameters with $v_1,\dots,v_n$ barred from $(-\infty,0] \cup (1,\infty)$. 
    Then, the series 

        \begin{equation}
        G_N\left(z,A,B\, \bigg{|}\, \begin{matrix} u_1^{-1}N & \dots & u_m^{-1}N \\
        v_1^{-1}N & \dots & v_n^{-1}N \end{matrix}\right) = 1+\sum_{d=1}^\infty z^d \sum_{\lambda \in 
        \mathrm{Y}_N^d}
        s_\lambda(A)s_\lambda(B)\frac{(u_1^{-1}N)_\lambda \dots (u_m^{-1}N)_\lambda}{(v_1^{-1}N)_\lambda \dots (v_n^{-1}N)_\lambda}
        \end{equation}

    \noindent
    is well-defined in every rank $N \in \N$. Its coefficients are

        \begin{equation}
            \begin{split}
                \frac{(u_1^{-1}N)_\lambda \dots (u_m^{-1}N)_\lambda}{(v_1^{-1}N)_\lambda \dots (v_n^{-1}N)_\lambda}
                &= \prod_{\Box \in \lambda} 
                \frac{(u_1^{-1}N+c(\Box)) \dots (u_m^{-1}N+c(\Box))}{(v_1^{-1}N+c(\Box)) \dots (v_n^{-1}N+c(\Box))} \\
                &=\left(N^{m-n} \frac{v_1\dots v_n}{u_1 \dots u_m} \right)^d
                \prod_{\Box \in \lambda} 
                \frac{(1+\frac{u_1c(\Box)}{N}) \dots (1+\frac{u_mc(\Box)}{N})}{(1+\frac{v_1c(\Box)}{N}) \dots (1+\frac{v_nc(\Box)}{N})},
            \end{split}
        \end{equation}

    \noindent
    where $d$ is the number of cells in $\lambda$. Absorbing the factor outside the product into 
    $z$, we obtain a new series 

        \begin{equation}
        \label{eqn:RenormalizedHypergeometric}
            K_N\left(q,A,B\, \bigg{|}\, \begin{matrix} u_1 & \dots & u_m \\
        v_1 & \dots & v_n \end{matrix}\right)
        = 1+\sum_{d=1}^\infty q^d \sum_{\lambda \in \Y_N^d}s_\lambda(A)s_\lambda(B)
        \frac{[\frac{u_1}{N}]_\lambda \dots [\frac{u_m}{N}]_\lambda}{[\frac{v_1}{N}]_\lambda \dots [\frac{v_n}{N}]_\lambda}
        \end{equation}

    \noindent
    whose coefficients are determined by the renormalized content polynomial 

        \begin{equation}
            [x]_\lambda = \prod_{\Box \in \lambda} (1+xc(\Box)).
        \end{equation}

    \noindent
    This series is absolutely convergent on the domain 

        \begin{equation}
            \left\{(q,A,B) \in \C^{1+2N} \colon |q|\|A\|\|B\| < \left| \frac{u_1 \dots u_m}{v_1 \dots v_n} \right|\right\}            
        \end{equation}

    \noindent
    when $m=n$, and absolutely convergent on all of $\C^{1+2N}$ when $m<n$.

    \subsection{Hypergeometric moments and cumulants}
    According to the above, provided $m \leq n$, we have a 
    sequence 

        \begin{equation}
            K_N\left(q,U_N,U_N^*\, \bigg{|}\, \begin{matrix} u_1 & \dots & u_m \\
        v_1 & \dots & v_n \end{matrix}\right), \quad N \in \N,
        \end{equation}

    \noindent
    of random univariate functions defined and analytic in an open neighborhood of $q=0$ in $\C$, 
    with power series expansion 

        \begin{equation}
            K_N\left(q,U_N,U_N^*\, \bigg{|}\, \begin{matrix} u_1 & \dots & u_m \\
        v_1 & \dots & v_n \end{matrix}\right) = 1 + \sum_{d=1}^\infty q^d
        \sum_{\lambda \in \Y_N^d} |s_\lambda(U_N)|^2 
          \frac{[\frac{u_1}{N}]_\lambda \dots [\frac{u_m}{N}]_\lambda}{[\frac{v_1}{N}]_\lambda \dots [\frac{v_n}{N}]_\lambda}.
        \end{equation}

    \noindent
    Let

        \begin{equation}
            K_N\left(q\, \bigg{|}\, \begin{matrix} u_1 & \dots & u_m \\
        v_1 & \dots & v_n \end{matrix}\right)= \mathbb{E} K_N\left(q,U_N,U_N^*\, \bigg{|}\, \begin{matrix} u_1 & \dots & u_m \\
        v_1 & \dots & v_n \end{matrix}\right), \quad N \in \N,
        \end{equation}

    \noindent
    be the corresponding sequence of expected values, a sequence of deterministic analytic functions, and let 

        \begin{equation}
              L_N\left(q\, \bigg{|}\, \begin{matrix} u_1 & \dots & u_m \\
        v_1 & \dots & v_n \end{matrix}\right) = \log K_N\left(q\, \bigg{|}\, \begin{matrix} u_1 & \dots & u_m \\
        v_1 & \dots & v_n \end{matrix}\right)
        \end{equation}

    \noindent
    be the principal logarithm of the expected value. 
    Write 

        \begin{equation}
             K_N\left(q\, \bigg{|}\, \begin{matrix} u_1 & \dots & u_m \\
        v_1 & \dots & v_n \end{matrix}\right)= 1+ \sum_{d=1}^\infty 
        \frac{q^d}{d!}  K_N^d\begin{pmatrix} u_1 & \dots & u_m \\
        v_1 & \dots & v_n \end{pmatrix}
        \end{equation}

    \noindent
    and 

         \begin{equation}
             L_N\left(q\, \bigg{|}\, \begin{matrix} u_1 & \dots & u_m \\
        v_1 & \dots & v_n \end{matrix}\right)= \sum_{d=1}^\infty 
        \frac{q^d}{d!}  L_N^d\begin{pmatrix} u_1 & \dots & u_m \\
        v_1 & \dots & v_n \end{pmatrix}
        \end{equation}

    \noindent
    for the Maclaurin expansions of these deterministic analytic functions. 
    We are going to derive $N \to \infty$ asymptotic expansions for the 
    hypergeometric moments

        \begin{equation}
             K_N^d\begin{pmatrix} u_1 & \dots & u_m \\
                v_1 & \dots & v_n \end{pmatrix}
        \end{equation}

    \noindent
    and hypergeometric cumulants 

        \begin{equation}
            L_N^d\begin{pmatrix} u_1 & \dots & u_m \\
        v_1 & \dots & v_n \end{pmatrix}
        \end{equation}

    \noindent
    in terms of disconnected and connected monotone simple Hurwitz numbers 
    with torus target. 

    \subsection{Asymptotic expansion of moments}
    The case of two upper and two lower parameters, which corresponds to the multivariate
    Gaussian hypergeometric function, exhibits all the features of the 
    general case. In the random matrix scaling, the multivariate hypergeometric function

        \begin{equation}
        G_N\left(z,A,B\, \bigg{|}\, \begin{matrix} x_1 & x_2 \\
        y_1 & y_2 \end{matrix}\right) = 1+\sum_{d=1}^\infty z^d \sum_{\lambda \in 
        \mathrm{Y}_N^d}
        s_\lambda(A)s_\lambda(B)\frac{(x_1)_\lambda (x_2)_\lambda}{(y_1)_\lambda  (y_2)_\lambda},
        \end{equation}

    \noindent
    renormalizes to

        \begin{equation}
        K_N\left(q,A,B\, \bigg{|}\, \begin{matrix} u_1 & u_2 \\
        v_1 & v_2 \end{matrix}\right) = 1+\sum_{d=1}^\infty q^d \sum_{\lambda \in 
        \mathrm{Y}_N^d} s_\lambda(A)s_\lambda(B)  
        \frac{[\frac{u_1}{N}]_\lambda [\frac{u_2}{N}]_\lambda}{[\frac{v_1}{N}]_\lambda[\frac{v_2}{N}]_\lambda}.
        \end{equation}

    \noindent
    This gives the random function

          \begin{equation}
          \label{eqn:RandomExample}
        K_N\left(q,U_N,U_N^*\, \bigg{|}\, \begin{matrix} u_1 & u_2 \\
        v_1 & v_2 \end{matrix}\right) = 1+\sum_{d=1}^\infty q^d \sum_{\lambda \in 
        \mathrm{Y}_N^d} |s_\lambda(U_N)|^2  
        \frac{[\frac{u_1}{N}]_\lambda [\frac{u_2}{N}]_\lambda}{[\frac{v_1}{N}]_\lambda[\frac{v_2}{N}]_\lambda},
        \end{equation}

    \noindent
    which is defined and analytic for $|q| < |\frac{u_1u_2}{v_1v_2}|$.
    By character orthogonality

        \begin{equation}
            \E|s_\lambda(U_N)|^2=1,
        \end{equation}

    \noindent
    hence by analyticity the expected value of \eqref{eqn:RandomExample} is

         \begin{equation}
        K_N\left(q \bigg{|}\, \begin{matrix} u_1 & u_2 \\
        v_1 & v_2 \end{matrix}\right) = 1+\sum_{d=1}^\infty q^d \sum_{\lambda \in 
        \mathrm{Y}_N^d} 
        \frac{[\frac{u_1}{N}]_\lambda [\frac{u_2}{N}]_\lambda}{[\frac{v_1}{N}]_\lambda[\frac{v_2}{N}]_\lambda}.
        \end{equation}

    \noindent
    We now analyze the $N \to \infty$ behavior of the finite sum 

        \begin{equation}
             K_N^d\left(\begin{matrix} u_1 & u_2 \\
            v_1 & v_2 \end{matrix}\right)=
             d!\sum_{\lambda \in \mathrm{Y}_N^d}
             \frac{[\frac{u_1}{N}]_\lambda [\frac{u_2}{N}]_\lambda}{[\frac{v_1}{N}]_\lambda[\frac{v_2}{N}]_\lambda}
        \end{equation}

    \noindent
    with $d \in \N$ arbitrary but fixed.

    The renormalized content polynomial of a Young diagram is given by

        \begin{equation}
            [x]_\lambda = \sum_{r=0}^{d-1}  x^r e_r(\lambda),
        \end{equation}

    \noindent
    where $e_r(\lambda) = e_r(c(\Box) \colon \Box \in \lambda)$ denotes
    the elementary symmetric polynomial of degree $r$ in $d$ variables
    evaluated on the multiset of contents of a Young diagram $\lambda$
    with $d$ cells. The sum terminates at $d-1$ because every Young diagram
    has a cell of content zero. Similarly, for any Young diagram $\lambda$
    with $d$ cells we have

        \begin{equation}
            \frac{1}{[x]_\lambda}  = \sum_{r=0}^\infty (-x)^r h_r(\lambda),           
        \end{equation}

    \noindent
    where $h_r(\lambda) = h_r(c(\Box) \colon \Box \in \lambda)$ is the complete
    symmetric polynomial of degree $r$ in $d$ variables evaluated on the content
    multiset of $\lambda$, and the series is absolutely convergent provided 
    $|x|c(\Box)<1$ for all $\Box \in \lambda$, in particular for $|x| < \frac{1}{d-1}$.
    Therefore,

        \begin{equation}
             K_N^d\left(\begin{matrix} u_1 & u_2 \\
            v_1 & v_2 \end{matrix}\right)= \sum_{r=0}^\infty \frac{1}{N^r}
            \sum_{\substack{(s_1,s_2,t_1,t_2) \in \N_0^4 \\ \|(s_1,s_2,t_1,t_2)\|_1=r}}
             u_1^{s_1}u_2^{s_2}(-v_1)^{t_1}(-v_2)^{t_2}d!
             \sum_{\lambda \in \Y_N^d} f_{(s_1,s_2,t_1,t_2)}(\lambda),
        \end{equation}

    \noindent
    where 

        \begin{equation}
                    f_{(s_1,s_2,t_1,t_2)}(\lambda) = e_{s_1}(\lambda)e_{s_2}(\lambda)h_{t_1}(\lambda)h_{t_2}(\lambda),
        \end{equation}

    \noindent
    and the series converges for all $N >(d-1)\max\{|v_1|,|v_2|\}$.

    Now consider the innermost sum above, which for $N \geq d$ is over the full set $\Y^d$ of
    Young diagrams with $d$ cells, which parameterize the irreducible representations $\V^\lambda$
    of the symmetric group $\group{S}^d$. For each $\lambda \in \Y^d$, the number 
    $f_{(s_1,s_2,t_1,t_2)}(\lambda)$ is the eigenvalue of a central element in the 
    group algebra $\C\group{S}^d$ acting in $\V^\lambda$. Explicitly, this element is

        \begin{equation}
            e_{s_1}(J_1,\dots,J_d)e_{s_2}(J_1,\dots,J_d)h_{t_1}(J_1,\dots,J_d)h_{t_2}(J_1,\dots,J_d),
        \end{equation}

    \noindent
    where 

        \begin{equation}
            \begin{split}
                J_1 &=0 \\
                J_2 &= (1\ 2) \\
                J_3 &= (1\ 3) + (2\ 3) \\
                \vdots \\
                J_d &= (1\ d) + \dots + (d-1\ d)
            \end{split}
            \end{equation}

    \noindent
    are the Jucys-Murphy elements in the group algebra \cite{DG,MN,OV}. Moreover,
    the commutator sum 

        \begin{equation}
            C = \sum_{\pi_1,\pi_2 \in \group{S}^d} \pi_1\pi_2\pi_1^{-1}\pi_2^{-1}
        \end{equation}

    \noindent
    is central and acts in $\V^\lambda$ as a scalar operator with eigenvalue $C(\lambda) = 
    (\frac{\dim \V^\lambda}{d!})^2$. Thus,

        \begin{equation}
            d!\sum_{\lambda \in \Y^d} f_{(s_1,s_2,t_1,t_2)}(\lambda)=
            \sum_{\lambda \in \Y^d}\frac{(\dim \V^\lambda)^2}{d!} C(\lambda)f_{(s_1,s_2,t_1,t_2)}(\lambda)
        \end{equation}

    \noindent
    is the normalized character of the central element 

        \begin{equation}
            C e_{s_1}(J_1,\dots,J_d)e_{s_2}(J_1,\dots,J_d)h_{t_1}(J_1,\dots,J_d)h_{t_2}(J_1,\dots,J_d)
        \end{equation}

    \noindent
    acting in the regular representation of $\C\group{S}^d$, which is the coefficient of 
    $\iota$ in the expansion of this product as a linear combination of conjugacy classes. 
    This is the number of factorizations of the form 

        \begin{equation}
        \label{eqn:MonotonePattern}
        \begin{split}
            \iota =& \pi_1\pi_2\pi_1^{-1}\pi_2^{-1} \\
            & (j_1^{(1)}\ k_1^{(1)}) \dots (j_{s_1}^{(1)}\ k_{s_1}^{(1)}), \quad k_1^{(1)} < \dots < k_{s_1}^{(1)} \\
            & (j_1^{(2)}\ k_1^{(2)}) \dots (j_{s_2}^{(2)}\ k_{s_2}^{(2)}), \quad k_1^{(2)} < \dots < k_{s_2}^{(2)} \\
            & (j_1^{(3)}\ k_1^{(3)}) \dots (j_{t_1}^{(3)}\ k_{t_1}^{(3)}), \quad k_1^{(3)} \leq \dots \leq k_{t_1}^{(3)} \\
            & (j_1^{(4)}\ k_1^{(4)}) \dots (j_{t_2}^{(4)}\ k_{t_2}^{(4)}), \quad k_1^{(4)} \leq \dots \leq k_{t_2}^{(4)},
        \end{split}
        \end{equation}

    \noindent
    where all transpositions $(j\ k)$ are in standard form $j<k$. This is bounded by
    the number of factorizations of $\iota$ of the form $\pi_1\pi_2\pi_1^{-1}\pi_2^{-1}\tau_1 \dots \tau_r$
    with $\tau_1,\dots,\tau_r$ any transpositions, and these unrestricted factorizations count possibly disconnected branched covers of a torus 
    with $r$ simple branch points. Consequently, a factorization of the form 
    \eqref{eqn:MonotonePattern} can only exist if $r=s_1+s_2+t_1+t_2$ satisfies
    the Riemann-Hurwitz constraint that $r=2g-2$ for $g \in \N$. In this case
    we denote the number of factorizations by $\mon_g^{d\bullet}(s_1,s_2;t_1,t_2)$, where the 
    bullet stands for disconnected. The argument is exactly the same in the general case.
    
    \begin{thm}
    \label{thm:AsymptoticDisconnected}
        For each $d \in \N$, we have 

            \begin{equation*}
                 K_N^d\begin{pmatrix} u_1 & \dots & u_m \\
        v_1 & \dots & v_n \end{pmatrix} \sim \sum_{g=1}^\infty N^{2-2g} 
        \mon_g^{d\bullet}\begin{pmatrix} u_1 & \dots & u_m \\
        v_1 & \dots & v_n \end{pmatrix}
            \end{equation*}

        \noindent
        as $N \to \infty$, where 

            \begin{equation*}
                \mon_g^{d\bullet}\begin{pmatrix} u_1 & \dots & u_m \\
        v_1 & \dots & v_n \end{pmatrix} = 
        \sum_{\substack{(s_1,\dots,s_m,t_1,\dots,t_n) \in \N_0^{m+n} \\ \|(s_1,\dots,s_m,t_1,\dots,t_n)\|_1=2g-2}} 
                 u_1^{s_1}\dots u_m^{s_m}
                 (-v_1)^{t_1} \dots (-v_n)^{t_n}
        \mon_g^{d\bullet}(s_1,\dots,s_m;t_1,\dots,t_n)
            \end{equation*}

        \noindent
        and $ \mon_g^{d\bullet}(s_1,\dots,s_m;t_1,\dots,t_n)$ counts factorizations
        of the form

            \begin{equation*}
        \begin{split}
            \iota =& \pi_1\pi_2\pi_1^{-1}\pi_2^{-1} \\
            & (j_1^{(1)}\ k_1^{(1)}) \dots (j_{s_1}^{(1)}\ k_{s_1}^{(1)}), \quad k_1^{(1)} < \dots < k_{s_1}^{(1)} \\
                &\vdots \\
            & (j_1^{(m)}\ k_1^{(m)}) \dots (j_{s_m}^{(m)}\ k_{s_m}^{(m)}), \quad k_1^{(m)} < \dots < k_{s_m}^{(m)} \\
            & (j_1^{(m+1)}\ k_1^{(m+1)}) \dots (j_{t_1}^{(m+1)}\ k_{t_1}^{(m+1)}), \quad k_1^{(m+1)} \leq \dots \leq k_{t_1}^{(m+1)} \\
                \vdots \\
            & (j_1^{(m+n)}\ k_1^{(m+n)}) \dots (j_{t_n}^{(m+n)}\ k_{t_n}^{(m+n)}), \quad k_1^{(m+n)} \leq \dots \leq k_{t_n}^{(m+n)},
        \end{split}
        \end{equation*}

    \end{thm}

    \subsection{Asymptotic expansion of cumulants}
    By a nonstandard application of the Exponential Formula explained in the next section,
    the counterpart of Theorem \ref{thm:AsymptoticDisconnected} for hypergeometric cumulants is
    the following. 

        \begin{thm}
        \label{eqn:AsymptoticConnected}
       For each $d \in \N$, we have

            \begin{equation*}
                 L_N^d\begin{pmatrix} u_1 & \dots & u_m \\
        v_1 & \dots & v_n \end{pmatrix} \sim \sum_{g=1}^\infty N^{2-2g} 
        \mon_g^{d}\begin{pmatrix} u_1 & \dots & u_m \\
        v_1 & \dots & v_n \end{pmatrix}
            \end{equation*}

        \noindent
        as $N \to \infty$, where 

            \begin{equation*}
                \mon_g^{d}\begin{pmatrix} u_1 & \dots & u_m \\
        v_1 & \dots & v_n \end{pmatrix} = 
        \sum_{\substack{(s_1,\dots,s_m,t_1,\dots,t_n) \in \N_0^{m+n} \\ \|(s_1,\dots,s_m,t_1,\dots,t_n)\|_1=2g-2}} 
                 u_1^{s_1}\dots u_m^{s_m}
                 (-v_1)^{t_1} \dots (-v_n)^{t_n}
        \mon_g^{d}(s_1,\dots,s_m;t_1,\dots,t_m)
            \end{equation*}

        \noindent
        and $\mon_g^{d}(s_1,\dots,s_m;t_1,\dots,t_m)$ counts the same factorizations
        as $\mon_g^{d\bullet}(s_1,\dots,s_m;t_1,\dots,t_m)$ with the additional constraint
        that the factors generate a transitive subgroup of $\group{S}^d$.
        \end{thm}

\section{Quasimodularity}
In this section we prove the general version of Theorem \ref{thm:QuasimodularCommutator}.

    \subsection{Bloch-Okounkov theorem}
    Let $\Lambda$ be the algebra of symmetric functions over $\C$, 
    and let $\Y$ be the set of all Young diagrams together with an
    empty diagram $\epsilon$. Write $|\lambda|$ for the number of cells
    in nonempty diagram and set $|\epsilon|=1$.
    Every $f \in \Lambda$ gives a corresponding function $f \colon \Y \to \C$ 
    on the set of Young diagrams defined by 

        \begin{equation}
            f(\lambda) = f(c(\Box) \colon \Box \in \lambda)
        \end{equation}

    \noindent
    for $\lambda$ nonempty and $f(\epsilon)=1$.
    The set of all such functions together with the zero function forms a unital subalgebra
    $\mathcal{C}(\Y)$ of the algebra of all functions $\Y \to \C$. The remarkable Bloch-Okounkov theorem 
    \cite{BlOk,Zagier} says that the $q$-average

        \begin{equation}
            \langle f \rangle_q = 
            \frac{\sum_{\lambda \in \Y} q^|\lambda f(\lambda)}{\sum_{\lambda \in \Y} q^{|\lambda|}}
        \end{equation}

    \noindent
    of $f \in \mathcal{C}(\Y)$
    is a polynomial in $E_2(q),E_4(q),E_6(q)$, i.e. a quasimodular form.
    In fact, the Bloch-Okounkov theorem holds for a larger algebra of functions on $\Y$,
    but the version stated above is sufficient for our purposes. We will use this
    result to prove that a version of Theorem \ref{thm:QuasimodularCommutator}
    holds for all hypergeometric cumulants of the CUE. We essentially follow 
    the argument in \cite{HIL}, which establishes a different generalization of 
    Theorem \ref{thm:QuasimodularCommutator}.
    
    \subsection{Total generating functions}
    Let 

        \begin{equation}
            \phi(q) = \frac{1}{\sum_{\lambda \in \Y} q^{|\lambda|}} = \prod_{n=1}^\infty (1-q^n).
        \end{equation}

    For each $g \geq 1$, let 

        \begin{equation}
            \vec{E}_g(q) = \sum_{d=1}^\infty \frac{q^d}{d!}
         \mon_g^{d\bullet}\begin{pmatrix} u_1 & \dots & u_m \\
        v_1 & \dots & v_n \end{pmatrix}
        \end{equation}

    \noindent
    be the generating function for the $g$th order contribution to the 
    asymptotic expansion of hypergeometric moments given in Theorem \ref{thm:AsymptoticDisconnected}.
    We omit the dependence of $\vec{E}_g(q)$ on the parameters $u_1,\dots,u_m,v_1,\dots,v_n$
    in order to lighten the notation. Note that

        \begin{equation}
            \vec{E}_1(q) = \prod_{n=1}^\infty \frac{1}{1-q^n},
        \end{equation}

    \noindent
    so that $\phi(q)\vec{E}_1(q)=1$.
    For $g \geq 2$, we have that $\phi(q)\vec{E}_g(q)$ is quasimodular
    by the Bloch-Okounkov theorem. 

    For each $g \geq 1$, let 

        \begin{equation}
            \vec{F}_g(q) = \sum_{d=1}^\infty \frac{q^d}{d!}
         \mon_g^{d}\begin{pmatrix} u_1 & \dots & u_m \\
        v_1 & \dots & v_n \end{pmatrix}
        \end{equation}

    \noindent
    be the generating function of the $g$th order contribution to the 
    asymptotic expansion of hypergeometric cumulants given in Theorem \eqref{eqn:AsymptoticConnected},
    and note that

        \begin{equation}
            \vec{F}_1(q) = \sum_{n=1}^\infty \log \frac{1}{1-q^n},
        \end{equation}

    \noindent
    so that $\vec{F}_1(q) = - \log \phi(q)$.
    
    Let $\hbar$ be a formal variable representing
    $1/N$ with $N=\infty$, and form the total generating functions

        \begin{equation}
            \vec{E}(q,\hbar) = \sum_{g=1}^\infty \hbar^{2g-2} 
        \vec{E}_g(q)
        \end{equation}

    \noindent
    and 

        \begin{equation}
            \vec{F}(q,\hbar) = \sum_{g=1}^\infty \hbar^{2g-2}
        \vec{F}_g(q).
        \end{equation}

    \noindent
    which satisfy

        \begin{equation}
        \label{eqn:ExponentialFormula}
           \log\left(1+ \vec{E}(q,\hbar)\right) = \vec{F}(q,\hbar),
        \end{equation}

    \noindent
    in the ring of formal power series $\C[[q,\hbar]]$.
    Remark that the use of the Exponential Formula here is not automatic 
    since $\hbar$ is an ordinary (rather than exponential) marker for genus, and 
    the monotonicity of the Hurwitz numbers in play is crucial; see \cite{GGN1,GGN2}.
    Indeed, this is fundamental to the relationship between unitary matrix integrals
    and monotone Hurwitz numbers, as the the $1/N$ expansion is an ordinary as opposed
    to exponential series in this perturbative parameter. 

    Now consider the series

        \begin{equation}
            \phi(q)\left(1+ \vec{E}(q,\hbar)\right) = 1 + \phi(q) + \sum_{g=2}^\infty \hbar^{2g-2}
            \phi(q) \vec{E}_g(q).
        \end{equation}

    \noindent
    This series has a well-defined logarithm in $\C[[\hbar,q]]$ such that 
    the coefficient of every positive power of $\hbar$ is quasimodular. 
    On the other hand, we have

        \begin{equation}
            \log \phi(q)(1+\vec{E}(q,\hbar)) = \log \phi(q) + \vec{F}(q,\hbar)
            = \sum_{g \geq 2} \hbar^{2g-2}\vec{F}_g(q),
        \end{equation}

    \noindent
    showing that $\vec{F}_g(q)$ is quasimodular for each $g \geq 2$.

\end{document}